\documentclass[11pt]{amsart}
\usepackage{amssymb,amsmath,amscd,amsfonts,amsthm,mathrsfs}
\usepackage{verbatim,paralist}
\usepackage{enumerate}
\RequirePackage{ifthen}
\usepackage[english]{babel}
\usepackage{xcolor}
 \usepackage{tikz}
\usepackage{latexsym}
\usepackage{amsmath}
\usepackage{amsfonts}
\usepackage{amssymb}
\usepackage{latexsym}

\usepackage{color}

\newcommand{\Qc}{\mathcal{Q}}

\newcommand{\Bc}{\mathcal{B}}

\newcommand{\Oc}{\mathcal{O}}

\newtheorem{theorem}{Theorem}[section]
\newtheorem{proposition}[theorem]{Proposition}
\newtheorem{lemma}[theorem]{Lemma}
\newtheorem{remark}[theorem]{Remark}
\newtheorem{example}[theorem]{Example}
\newtheorem{definition}[theorem]{Definition}
\newtheorem{corollary}[theorem]{Corollary}

\title{ Riesz projection and essential $S$-spectrum in quaternionic setting}
\author{Hatem Baloudi}
\address{Hatem Baloudi, Department of Mathematics, Faculty of Sciences of Gafsa, University of Gafsa, 2112 Zarroug,
Tunisia.\\
E-mail adress: hatem.beloudi@gmail.com; hatem.beloudi@fsg.u-gafsa.tn}
\author{Sayda Belgacem}
\address{Sayda Belgacem, Departement of Mathematics Faculty of Sciences of Sfax, University of Sfax, Route de Soukra km 3.5, B.P. 1171, 3000 Sfax, Tunisia.\\
 E-mail adress: saidabelgacem207@yahoo.fr}
\author{Aref Jeribi}
\address{Aref Jeribi, Departement of Mathematics Faculty of Sciences of Sfax, University of Sfax, Route de Soukra km 3.5, B.P. 1171, 3000 Sfax, Tunisia.\\
aref.jeribi@gmail.com}


\begin{document}

\maketitle

\begin{abstract}
This paper is devoted to the investigation of the Weyl and the essential $S-$spectra of a bounded right quaternionic linear operator in a right quaternionic Hilbert space. Using the quaternionic Riesz projection, the $S-$eigenvalue of finite type is both introduced and studied. In particular, we have shown that the Weyl and the essential $S-$spectra do not contain eigenvalues of finite type. We have also described the boundary of the Weyl $S-$spectrum and the particular case of the spectral theorem of the essential $S-$spectrum.
\end{abstract}
\tableofcontents
\vskip 0.2 cm

\section{Introduction}
\noindent Over the recent years, the spectral theory for quaternionic operators has piqued the interest and attracted the attention of multiple researchers, see for instance \cite{SL,KRE,NCFCBO,FJDP,FJGanter,2021,FIF,FIDC,BK} and references therein. Research in this topic is motivated by application in various fields, including quantum mechanics, fractional evolution problems \cite{FJGanter}, and quaternionic Schur analysis \cite{DFIS}. The concept of spectrum is one of the main objectives in the theory of quaternionic operators acting on quaternionic Hilbert spaces. The obscurity that seemed to surround the precise definition of the quaternionic spectrum of a linear operator was deemed a stumbling block. However, in 2006, F. Colombo and I. Sabadini succeeded in devising a new notion conducive to the complete development of the quaternionic operator theory, namely the $S$-spectrum. We refer to \cite[Subsection 1.2.1]{FJDP} for a precise history. After several more years of formulating the spectral theorem that stemmed from the $S$-spectrum, D. Alpay, F. Colombo, D.P. Kimsey, see \cite{DFDP}, supplied ample evidence in 2016 to further establish this fundamental theorem for both bounded and unbounded operators. In the book \cite{FJDP}, see also \cite{FJGanter}, the authors briefly explain the concept of $S-$spectrum and give the systematic basis of quaternionic spectral theory. We refer to the book \cite{DFIS} and the references therein for the spectral theory on the $S$-spectrum for Clifford operators and to \cite{PF} for some results on operators perturbation.
\vskip 0.2 cm

Motivated by the new concept of $S-$spectrum in the quaternionic setting, Muraleetharan and Thirulogasanthar, in \cite{BK,BK2}, introduced the Weyl and essential $S$-spectra and gave a characterization using Fredholm operators. We refer to \cite{B1} for the study of the general framework of the Fredholm element with respect to a quaternionic Banach algebra homomorphism. In general, the set of all operators acting on right Banach space is not quaternionic Banach algebra with respect to the composition operators. By \cite[Theoreme 7.1 and Theorem 7.3]{RVA}, if $V_{\mathbb{H}}^{R}$ is a separable quaternionic Hilbert space, then $\mathcal{B}(V_{\mathbb{H}}^{R})$ (the set of all right bounded operators) is a quaternionic two-sided $C^{*}-$algebra and the set of all compact operators $\mathcal{K}(V_{\mathbb{H}}^{R})$ is a closed two-sided ideal of $\mathcal{B}(V_{\mathbb{H}}^{R})$ which is closed under adjunction. In this regard, in \cite{BK}, the author defined the essential $S-$spectrum as the $S-$spectrum of quotient map image of bounded right linear operator on the Calkin algebra $\mathcal{B}(V_{\mathbb{H}}^{R})/\mathcal{K}(V_{\mathbb{H}}^{R})$.
\vskip 0.2 cm

In order to explain the objective of this work, we start by recalling a few results concerning the discrete spectrum and the Riesz projection in the complex setting. Let $T$ be a linear operator acting on a complex Banach space $V_{\mathbb{C}}$. We denote the spectrum of $T$ by $\sigma(T)$. Let $\sigma$ be an isolated part of $\sigma(T)$. The Riesz projection of $T$ corresponding to $\sigma$ is the operator
\begin{align*}P_{\sigma}=\displaystyle\frac{1}{2\pi i}\int_{C_{\sigma}}(z-T)^{-1}dz\end{align*}
where $C_{\sigma}$ is a smooth closed curves belonging  to the resolvent set $\mathbb{C}\backslash\sigma(T)$ such that $C_{\sigma}$ surrounds $\sigma$ and separates $\sigma$ from $\sigma(T)\backslash\sigma$. The discrete spectrum of $T$, denoted $\sigma_{d}(T)$, is the set of isolated point $\lambda\in\mathbb{C}$ of $\sigma(T)$ such that the corresponding Riesz projection $P_{\{\lambda\}}$ are finite dimensional, see \cite{Gohberg,Lutgen}. Note that in general we have $\sigma_{e}(T)\subset\sigma(T)\backslash\sigma_{d}(T)$, where $\sigma_{e}(T)$ denotes the set of essential spectrum of $T$. We refer to \cite{BJ,J1,Wolf} for more properties of $\sigma_{e}(T)$. Note that, if $A$ is a self-adjoint operator on a Hilbert space, then $\sigma_{e}(T)=\sigma(T)\backslash\sigma_{d}(T)$. In particular, the essential spectrum is empty if and only if $\sigma(T)=\sigma_{d}(T)$. We point out that this point, namely the absence of the essential spectrum, has been studied in many works, e.g., \cite{Gol,Keller}.
\vskip 0.2 cm

In the quaternionic setting, if $T\in\mathcal{B}(V_{\mathbb{H}}^{R})$ and $q\in\sigma_{S}(T)\backslash\mathbb{R}$ (where $\sigma_{S}(T)$ denote the $S-$spectrum of $T$), then $q$ is not an isolated point of $\sigma_{S}(T)$. Indeed, $[q]:=\{hqh^{-1}:\ h\in\mathbb{H}^{*}\}\subset \sigma_{S}(T)$, see \cite{FJDP}. By the compactness of $[q]$, $\{q\}$ is not an isolated part of $\sigma_{S}(T)$. However, if we set
\begin{align*}\Omega:=\sigma_{S}(T)/\sim\end{align*}
where $p\sim q$ if and only if $p\in [q]$, and $E_{T}$ the set of representative, then if $[q]$ is an isolated part of $\sigma_{S}(T)$ then $q$ is an isolated point of $E_{T}$.
\vskip 0.2 cm

The first aim of this work is to study the isolated part of the $S-$ spectrum of a bounded right quaternionic operators and its relation with the essential $S$-spectrum. To begin with, we consider the Riesz projector associated with a given quaternionic operator $T$ which was introduced in \cite{FIDC1}. We refer to \cite{DFJI,KRE,PF,FJDP} for more details on this concept. We treat the decomposition of the essential $S-$spectrum of $T$ as a function of the Riesz projector, and more generally as a function of a projector which commutes with $T$, see Theorem \ref{t:r}. The technique of the proof is inspired from  \cite{KRE}. We also discuss the Riesz decomposition theorem \cite[Theorem 6]{PSK} in quaternionic setting. More precisely, we prove that this decomposition is unique. Motivated by this, we study the quaternionic version of the discrete $S-$spectrum. Following the complex formalism given in \cite{Gohberg, Lutgen}, we show that the essential $S-$spectrum of a given right operator acting on right Hilbert space does not contain discrete element of the $S-$spectrum. The second aim of this work is to give  new results concerning the Weyl and essential $S$-spectra in quaternionic setting. First off, Theorem \ref{t:4} gives a description of the boundary of the Weyl $S-$spectrum. The proof is based on the study of the minimal modulus  of the right quaternionic operators. We also deal with the particular case of the spectral theorem of essential  $S-$spectra. The technique of the proof is inspired from \cite{FIDC1}.
\vskip 0.1 cm

The article is organised as follows: In Section \ref{sec:1}, we present general definitions about operators theory in quaternionic setting. In Section \ref{sec:2}, we discuss the question of decomposition of the essential $S-$spectrum. Finally, in Section \ref{sec:3}, we provide  new results of the Weyl $S-$spectrum.

\section{Mathematical preliminaries}\label{sec:1}
In this section, we review some basic notions about quaternions, right quaternionic Hilbert space, right linear operator (even unbounded and define its $S-$spectrum), and slice functional calculus. For details, we refer to the reader  \cite{SL,FIDC1,FJDP,RVA,WR}.
\subsection{Quaternions}
Let $\mathbb{H}$ be the Hamiltonian skew field of quaternions. This class of numbers can be written as:
\begin{align*}
q=q_{0}+q_{1}i+q_{2}j+q_{3}k,
\end{align*}
where $q_{l}\in \mathbb{R}$ for $l=0,1,2,3$ and $i\ ,j\ ,k$ are the three quaternionic imaginary units satisfying
\begin{align*}
i^{2}=j^{2}=k^{2}=ijk=-1.
\end{align*}
The real and the imaginary part of $q$ is defined as ${\rm Re}(q)=q_{0}$ and ${\rm Im}(q)=q_{1}i+q_{2}j+q_{3}k$, respectively. Then, the conjugate and the usual norm of the quaternion q are given respectively by
\begin{align*}
\overline{q}=q_{0}-q_{1}i-q_{2}j-q_{3}k \mbox{ and  }|q|=\sqrt{q\overline{q}}.
\end{align*}
 The set of all imaginary unit quaternions in $\mathbb{H}$ is denoted by $\mathbb{S}$ and  defined as
 \begin{align*}
\mathbb{S}=\Big\{q_{1}i+q_{2}j+q_{3}k:\ q_{1},\ q_{2},\ q_{3}\in\mathbb{R}, q_{1}^{2}+q_{2}^{2}+q_{3}^{2}=1  \Big\}.
\end{align*}
The name imaginary unit is due the fact that, for any  $I\in\mathbb{S}$, we have
\begin{align*}I^{2}=-\overline{I}I=-|I|^{2}=-1.\end{align*}
\noindent  For every $q\in\mathbb{H}\backslash\mathbb{R}$, we associate the  unique element
\begin{align*}I_{q}:=\frac{{\rm Im}(q)}{|{\rm Im}(q)|}\in\mathbb{S}\end{align*}
 such that
\begin{align*}q={\rm Re}(q)+I_{q}|{\rm Im}(q)|.\end{align*}
This implies that
\begin{align*}\mathbb{H}=\bigcup_{I\in\mathbb{S}}\mathbb{C}_{I}.\end{align*}
where
\begin{align*}\mathbb{C}_{I}:=\mathbb{R}+I\mathbb{R}.\end{align*}
We can associate to $q\in\mathbb{H}$ the $2-$dimensional sphere
\begin{align*}[q]:=\{{\rm Re}(q)+I\vert {\rm Im}(q)\vert:\ I\in\mathbb{S}\}.\end{align*}
This sphere has center at the real point ${\rm Re}(q)$ and radius $\vert {\rm Im}(q)\vert$.
\subsection{Right quaternionic Hilbert space and operator}
In this subsection, we recall the concept of right quaternionic Hilbert space and right linear operator (see, \cite{SL,NCFCBO,RVA}).
\begin{definition}\cite{SL}
{\rm Let $V_{\mathbb{H}}^{R}$ be a right vector space. The map
\begin{align*}\langle .,.\rangle :V_{\mathbb{H}}^{R}\times V_{\mathbb{H}}^{R}\longrightarrow \mathbb{H}\end{align*}
is called an inner product if it satisfies the following properties:\\
\noindent $(i)$ $\langle f,gq+h\rangle=\langle f,g\rangle q+\langle f,h\rangle$, for all $f,g,h\in V_{\mathbb{H}}^{R}$ and $q\in\mathbb{H}$.\\
$(ii)$ $\langle f,g\rangle=\overline{\langle g,f\rangle}$, for all $f,g\in V_{\mathbb{H}}^{R}$.\\
$(iii)$ If $f\in V_{\mathbb{H}}^{R}$, then $\langle f,f\rangle\geq 0$ and $f=0$ if $\langle f,f\rangle=0$.\\
The pair $(V_{\mathbb{H}}^{R},\langle .,.\rangle)$ is called a right quaternionic pre-Hilbert space.  Moreover, $V_{\mathbb{H}}^{R}$ is said to be right quaternionic Hilbert space, if
\begin{align*}\|f\|=\sqrt{\langle f,f\rangle}\end{align*}
defines a norm for which $V_{\mathbb{H}}^{R}$ is complete.}
\end{definition}
\vskip 0.2 cm

In the sequel, we assume that $V_{\mathbb{H}}^{R}$ is complete and separable.
We now recall the concept of Hilbert basis in the quaternionic case. First, we review the following proposition, the proof of which is similar to its complex version, see \cite{RVA,KV}.
\begin{proposition}\label{p:2.2}
Let  $V_{\mathbb{H}}^{R}$ be a right quaternionic Hilbert space and  let $\mathcal{F}=\{f_{k}:\ k\in\mathbb{N}\}$ be an orthonormal subset of $V_{\mathbb{H}}^{R}$. The following properties are equivalent:
\begin{enumerate}
\item For every $f,g\in V_{\mathbb{H}}^{R}$, the series $\sum_{k\in\mathbb{N}}\langle f,f_{k}\rangle\langle f_{k},g\rangle$ converges absolutely and
\begin{align*} \langle f,g\rangle=\sum_{k\in\mathbb{N}}\langle f,f_{k}\rangle\langle f_{k},g\rangle.\end{align*}
\item For every $f\in V_{\mathbb{H}}^{R}$, we have
\begin{align*}\|f\|^{2}=\sum_{k\in\mathbb{N}}\vert \langle f_{k},f\rangle\vert^{2}\end{align*}
\item $\mathcal{F}^{\bot}:=\Big\{f\in V_{\mathbb{H}}^{R}:\ \langle f,g \rangle=0\mbox{ for all }g\in\mathcal{F}\Big\}=\{0\}$.
\item $\langle \mathcal{F}\rangle:=\Big\{\displaystyle\sum_{l=1}^{m}f_{l}q_{l}:\ f_{l}\in\mathcal{F},\ q_{l}\in\mathbb{H},\ m\in\mathbb{N}\Big\}$ is dense in $V_{\mathbb{H}}^{R}$.
\end{enumerate}
\end{proposition}
\begin{definition}
{\rm Let $\mathcal{F}$ be an orthonormal subset of $V_{\mathbb{H}}^{R}$. $\mathcal{F}$ is said to be Hilbert basis of $V_{\mathbb{H}}^{R}$ if $\mathcal{F}$ verifies one of the equivalent conditions of Proposition \ref{p:2.2}.}
\end{definition}
\vskip 0.2 cm

The proof of the following proposition is the same as its complex version, see \cite{RVA,KV}.
\begin{proposition}
Let  $V_{\mathbb{H}}^{R}$ be a right quaternionic Hilbert space. Then,
\begin{enumerate}
\item $V_{\mathbb{H}}^{R}$ admits a Hilbert basis.
\item Two Hilbert basis of $V_{\mathbb{H}}^{R}$ have the same cardinality.
\item If $\mathcal{F}$ is  a Hilbert basis of $V_{\mathbb{H}}^{R}$, then every $f\in V_{\mathbb{H}}^{R}$ can be uniquely decomposed as follows:
\begin{align*}f=\sum_{k\in\mathbb{N}}f_{k}\langle f_{k},f\rangle\end{align*}
where the series $\sum_{k\in\mathbb{N}}f_{k}\langle f_{k},f\rangle$ converges absolutely in $V_{\mathbb{H}}^{R}$.
\end{enumerate}
\end{proposition}
\vskip 0.2 cm

The quaternionic multiplication is not commutative. Afterwards, we recall that if $V_{\mathbb{H}}^{R}$ is a right separable quaternionic Hilbert space, we can define the left scalar multiplication on $V_{\mathbb{H}}^{R}$  using an arbitrary Hilbert basis on $V_{\mathbb{H}}^{R}$. We refer to \cite{RVA} for an explanation of this construction. Let $\mathcal{F}=\Big\{f_{k}:\ k\in\mathbb{N}\Big\}$ be a Hilbert basis of $V_{\mathbb{H}}^{R}$. The left scalar multiplication on $V_{\mathbb{H}}^{R}$ induced by $\mathcal{F}$ is defined as the map
 \begin{align*}
&\mathbb{H}\times V_{\mathbb{H}}^{R}\longrightarrow V_{\mathbb{H}}^{R}
\\
&\ (q,f)\longmapsto\ qf=\displaystyle\sum_{k\in\mathbb{N}}f_{k}q\langle f_{k},f\rangle.
\end{align*}
\vskip 0.2 cm

 The properties of the left scalar multiplication are described in the following proposition.
 \begin{proposition}\cite[Proposition 3.1]{RVA} Let $f,g\in V_{\mathbb{H}}^{R}$ and $p,q\in\mathbb{H}$, then
 \begin{enumerate}
 \item $q(f+g)=qf+qg$ and $q(fp)=(qf)p$.
 \item $\|qf\|=|q|\|f\|$.
 \item $q(pf)=(qp)f$.
 \item $\langle \overline{q}f,g\rangle=\langle f,qg\rangle$.
 \item $rf=fr$, for all $r\in\mathbb{R}$.
 \item $qf_{k}=f_{k}q$, for all $k\in\mathbb{N}$.
 \end{enumerate}
 \end{proposition}
 It is easy to see that $(p+q)f=pf+qf$, for all $p,q\in\mathbb{H}$ and $f\in V_{\mathbb{H}}^{R}$. In the sequel, we consider $V_{\mathbb{H}}^{R}$ as a right quaternionic Hilbert space equipped with the left scalar multiplication.
\begin{definition}{\rm Let $V_{\mathbb{H}}^{R}$ be a right quaternionic Hilbert space. A mapping $T:\mathcal{D}(T)\subset V_{\mathbb{H}}^{R}\longrightarrow V_{\mathbb{H}}^{R}$, where $\mathcal{D}(T)$ denote the domain of $T$, is called \emph{quaternionic right linear} if
\begin{align*}T(f+gq)=T(f)+T(g)q,\mbox{ for all }f,\ g\in\mathcal{D}(T) \mbox{ and }q\in\mathbb{H}.\end{align*}
The operator $T$ is called closed, if the graph $\mathcal{G}(T):=\{(f,Tf):\ f\in\mathcal{D}(T)\}$ is a closed right linear subspace of $V_{\mathbb{H}}^{R}\times V_{\mathbb{H}}^{R}$.}
\end{definition}
\noindent We call an quaternionic right operator $T$ bounded if
\begin{align*}\|T\|:=\sup\Big\{\|Tf\|:\ f\in V_{\mathbb{H}}^{R},\ \|f\|=1\Big\}<+\infty.\end{align*}
\vskip 0.2 cm

The set of all bounded right operators on $V_{\mathbb{H}}^{R}$ is denoted by $\mathcal{B}(V_{\mathbb{H}}^{R})$ and the identity operator on $V_{\mathbb{H}}^{R}$ will be denoted by $\mathbb{I}_{V_{\mathbb{H}}^{R}}$. Let $T\in \mathcal{B}(V_{\mathbb{H}}^{R})$, we denote the null space of $T$  by $N(T)$ and its  range space by $R(T)$. A closed subspace $M$ of $V_{\mathbb{H}}^{R}$ is said to be $T-$invariant subspace if $T(M)\subset M$. Note that the function $f\longmapsto Tf-fq$ is not right linear, we refer to \cite{NCFCBO, FJDP} for this point of view. The fundamental suggestion of \cite{NCFCBO} is to define the spectrum using the Cauchy kernel series. We recall these concepts from the book \cite{FJDP}.
\begin{definition}{\rm Let $T:\mathcal{D}(T)\subset V_{\mathbb{H}}^{R}\longrightarrow V_{\mathbb{H}}^{R}$ be a right linear operator. We define the operator $Q_{q}(T):\mathcal{D}(T^{2})\longrightarrow V_{\mathbb{H}}^{R}$ by
\begin{align*}Q_{q}(T):=T^{2}-2{\rm Re}(q)T+|q|^{2}\mathbb{I}_{V_{\mathbb{H}}^{R}}.\end{align*}
\noindent $1)$ The \emph{$S-$resolvent} set of $T$ is defined as follows:
\begin{align*}\rho_{S}(T):=\Big\{q\in\mathbb{H}:\ N(Q_{q}(T))=\{0\},\overline{R(T)}=\mathbb{H}\mbox{ and }Q_{q}(T)^{-1}\in\mathcal{B}(V_{\mathbb{H}}^{R})\Big\}.\end{align*}
\noindent $2)$ The \emph{$S-$spectrum} of $T$ is defined as:
\begin{align*}\sigma_{S}(T)=\mathbb{H}\backslash\rho_{S}(T).\end{align*}
\noindent $3)$ The \emph{point $S-$spectrum} of $T$ is given by
\begin{align*}\sigma_{pS}(T):=\Big\{q\in\mathbb{H}:\ N(Q_{q}(T))\neq\{0\}\Big\}.\end{align*}}
\end{definition}
\vskip 0.1 cm

 For $T\in\mathcal{B}(V_{\mathbb{H}}^{R})$, the $S-$spectrum $\sigma_{S}(T)$ is a non-empty compact set, see \cite{FJDP}. We recall that if $T\in\mathcal{B}(V_{\mathbb{H}}^{R})$ and $q\in\sigma_{S}(T)$, then all the elements of the sphere $[q]$ belong to $\sigma_{S}(T)$, see \cite[Theorem 7.2.8]{DFIS}.

  Let  $v\in V_{\mathbb{H}}^{R}\backslash\{0\}$, then $v$ is a right eigenvalue of $T$ if $T(v)$ is a right quaternionic multiple of $v$. That, is
\begin{align*}T(v)=vq\end{align*}
where $q\in\mathbb{H}$, known as the right eigenvalue. The set of right eigenvalues coincides with the point $S-$spectrum, see \cite[Proposition 4.5]{RVA}.
\vskip 0.1 cm

Now, we recall the definition of essential $S-$spectrum, we refer to \cite{BK,BK2} for more details. Let $V_{\mathbb{H}}^{R}$ be a separable right quaternionic Hilbert space equipped with a left scalar multiplication. Using \cite{RVA}, $\mathcal{B}(V_{\mathbb{H}}^{R})$ is a quaternionic two-sided Banach $C^{*}-$algebra with unity, and the set of all compact operators $\mathcal{K}(V_{\mathbb{H}}^{R})$ is a closed two-sided ideal of  $\mathcal{B}(V_{\mathbb{H}}^{R})$. We consider the natural quotient map:
 \begin{align*}
&\pi: \mathcal{B}(V_{\mathbb{H}}^{R}) \longrightarrow \mathcal{C}(V_{\mathbb{H}}^{R}):=\mathcal{B}(V_{\mathbb{H}}^{R})/\mathcal{K}(V_{\mathbb{H}}^{R})
\\
&\ \ \quad \ \quad T \longmapsto\ \ [T]=T+\mathcal{K}(V_{\mathbb{H}}^{R}).
\end{align*}
 Note that $\pi$ is a unital homomorphism, see \cite{BK}. The norm on $\mathcal{C}(V_{\mathbb{H}}^{R})$ is given by
 \begin{align*}\|[T]\|=\inf_{K\in\mathcal{K}(V_{\mathbb{H}}^{R})}\|A+K\|.\end{align*}
\begin{definition}\cite{BK}
{\rm The essential $S-$spectrum of $T\in\mathcal{B}(V_{\mathbb{H}}^{R})$ is the $S-$spectrum of $\pi(A)$ in the Calkin algebra $\mathcal{C}(V_{\mathbb{H}}^{R})$. That is,
\begin{align*}\sigma_{e}^{S}(T):=\sigma_{S}(\pi(A)).\end{align*}}
\end{definition}
\subsection{The quaternionic functional calculus}
The quaternionic functional calculus is defined on the class of slice regular function $f:U\longrightarrow\mathbb{H}$ for some set $U\subset \mathbb{H}$. We recall this concept and refer to \cite{PF,FIDC1,FJDP} and the references therein on the matter.
\begin{definition}
{\rm A set $U\subset\mathbb{H}$ is called\\
\noindent $(i)$ axially symmetric if $[x]\subset U$ for any $x\in U$ and\\
\noindent $(ii)$ a slice domain if $U$ is open, $U\cap\mathbb{R}\neq\emptyset$ and $U\cap\mathbb{C}_{I}$ is a domain in $\mathbb{C}_{I}$, for any $I\in\mathbb{S}$.}
\end{definition}
\begin{definition}{\rm  Let $U\subset\mathbb{H}$ be an open set.
A real differentiable function $f:U\longrightarrow \mathbb{H}$ is said to be left $s-$regular $($resp. right $s-$regular$)$ if for every $I\in\mathbb{S}$, the function $f$ satisfy
\begin{align*}\displaystyle\frac{1}{2}[\frac{\partial}{\partial x}f(x+Iy)+I\frac{\partial}{\partial y}f(x+Iy)]=0\ \Big(resp. \frac{1}{2}[\frac{\partial}{\partial x}f(x+Iy)\\
\quad\quad \quad\quad\quad+\frac{\partial}{\partial y}f(x+Iy)I]=0\Big).\end{align*}}
\end{definition}
\vskip 0.2 cm

We denote the class of left $s-$regular $($resp. right $s-$regular$)$ by $\mathcal{R}^{L}(U)$ $($resp. $\mathcal{R}^{R}(U))$. We recall that $\mathcal{R}^{L}(U)$ is a right $\mathbb{H}-$module and $\mathcal{R}^{R}(U)$ is a left $\mathbb{H}-$module. Let $V_{\mathbb{H}}^{R}$ be a separable right quaternionic Hilbert space equipped with a Hilbert basic $\mathcal{N}$ and with a left scalar multiplication. We recall that $\mathcal{B}(V_{\mathbb{H}}^{R})$ is a two-sided ideal quaternionic Banach algebra with respect to the left multiplication given by
\begin{align*}(q.T)f=\displaystyle \sum_{g\in\mathcal{N}}gq\langle g,Tf\rangle\mbox{ and }(Tq)f=\displaystyle\sum_{g\in\mathcal{N}}T(g)q\langle g,f\rangle.\end{align*}
\vskip 0.2 cm

\begin{definition}{\rm Let $T\in \mathcal{B}(V_{\mathbb{H}}^{R})$ and $q\in\rho_{S}(T)$. The \emph{left $S-$resolvent operator} is given by
\begin{align*}S_{L}^{-1}(q,T):=-Q_{q}(T)^{-1}(T-\overline{q}\mathbb{I}_{V_{\mathbb{H}}^{R}}),\end{align*}
and the \emph{right $S-$resolvent operator} is defined by
\begin{align*}S_{R}^{-1}(q,T):=-(T-\overline{q}\mathbb{I}_{V_{\mathbb{H}}^{R}})Q_{q}(T)^{-1}.\end{align*}}
\end{definition}
\begin{definition}\cite[Definition 3.4]{KRE}
{\rm Let $T\in \mathcal{B}(V_{\mathbb{H}}^{R})$ and let $U\subset\mathbb{H}$ be an axially symmetric $s-$domain that contains the $S-$spectrum $\sigma_{S}(T)$ and such that $\partial(U\cap\mathbb{C}_{I})$ is union of a finite number of continuously differentiable Jordan curves for every $I\in\mathbb{S}$. We say that $U$ is a $T-$admissible open set.}
\end{definition}
\begin{definition}
{\rm Let $T\in \mathcal{B}(V_{\mathbb{H}}^{R})$, $W\subset \mathbb{H}$ be open set. A function $f\in \mathcal{R}^{L}(W)($resp. $\mathcal{R}^{R}(W))$is said to be locally left regular $($resp. right regular$)$ function on $\sigma_{S}(T)$, if there is $T-$admissible domain $U\subset\mathbb{H}$ such that $\overline{U}\subset W$.}
\end{definition}
\vskip 0.2 cm

We  denote by $\mathcal{R}^{L}_{\sigma_{S}(T)}$ $($resp. $\mathcal{R}^{R}_{\sigma_{S}(T)})$ the set of locally left $($ resp. right$)$ regular functions on $\sigma_{S}(T)$. Now, we recall the two versions of the quaternionic functional calculus.
\begin{definition}\cite[Definition 4.10.4]{FIDC1}
{\rm Let $T\in \mathcal{B}(V_{\mathbb{H}}^{R})$ and $U\subset\mathbb{H}$ be a $T-$admissible domain. Then,
\begin{equation}\label{e:4}f(T)=\displaystyle\frac{1}{2\pi}\int_{\partial(U\cap\mathbb{C}_{I})}S_{L}^{-1}(q,T)dq_{I}f(q)\ \forall f\in \mathcal{R}^{L}_{\sigma_{S}(T)} \end{equation}
and
\begin{equation}\label{e:5}f(T)=\displaystyle\frac{1}{2\pi}\int_{\partial(U\cap\mathbb{C}_{I})}f(q)dq_{I}S_{R}^{-1}(q,T)\ \forall f\in \mathcal{R}^{R}_{\sigma_{S}(T)} \end{equation}
where $dq_{I}=-dqI$.}
\end{definition}
The two integrals that appear in Eqs $(\ref{e:4})$ and $(\ref{e:5})$ are independent of the choice of imaginary unit $I\in\mathbb{S}$ and $T-$admissible domain, see \cite[Theorem 4.10.3]{FIDC1}.
\vskip 0.1 cm

A set $\sigma$ is called an isolated part of $\sigma_{S}(T)$ if both $\sigma$ and $\sigma_{S}(T)\setminus\sigma$ are closed subsets of $\sigma_  {S}(T)$.
\begin{definition}\cite{A:C,PF}
{\rm Let $T\in \mathcal{B}(V_{\mathbb{H}}^{R})$. Denote by $U_{\sigma}$ an axially symmetric $s-$domain that contains the axially symmetric isolated part $\sigma\subset\sigma_{S}(T)$ but not any other point of $\sigma_{S}(T)$. Suppose that the Jordan curves $\partial(U_{\sigma}\cap\mathbb{C}_{I})$ belong to the $S-$resolvent set $\rho_{S}(T)$, for any $I\in\mathbb{S}$. We define the quaternionic Riesz projection by
\begin{align*}P_{\sigma}=\displaystyle\frac{1}{2\pi}\int_{\partial(U_{\sigma}\cap\mathbb{C}_{I})}S_{L}^{-1}(q,T)dq_{I}.\end{align*}}
\end{definition}
\begin{remark}\cite{A:C,PF}
{\rm The concept of $P_{\sigma}$ can be given by using the right $S-$resolvent operator $S_{R}^{-1}(q,T)$, that is
\begin{align*}P_{\sigma}=\displaystyle\frac{1}{2\pi}\int_{\partial(U_{\sigma}\cap\mathbb{C}_{I})}dq_{I}S_{R}^{-1}(q,T).\end{align*}}
\end{remark}
\vskip 0.1 cm Note that $P_{\sigma}$ is a projection  that  commutes with $T$, see \cite[Theorem 2.8]{A:C}.
\vskip 0.1 cm

In the sequel, we assume that $V_{\mathbb{H}}^{R}$ is a separable right quaternionic Hilbert space with infinite dimensional.
\section{Riesz projection and essential spectrum}\label{sec:2}
We recall that in \cite{BK,BK2}, the study of the essential $S-$spectrum is established using the theory of Fredholm operators. The aim of this section is to show that the essential $S-$spectrum does not contain discrete element of the $S-$spectrum. In this regard, let $V_{\mathbb{H}}^{R}$ be a separable right quaternionic Hilbert space. We note that if $K\in\mathcal{K}(V_{\mathbb{H}}^{R})$, then $\sigma_{e}^{S}(A+K)=\sigma_{e}^{S}(A)$ for all $A\in\mathcal{B}(V_{\mathbb{H}}^{R})$. We start by showing that in general the $S$-spectrum does not satisfy this property. We define $V_{\mathbb{H}}^{'}=\mathcal{B}(V_{\mathbb{H}}^{R},\mathbb{H})$ and call $V_{\mathbb{H}}^{'}$ the right dual space of $V_{\mathbb{H}}^{R}$.
\begin{theorem}\label{t:r}
Let $T\in \mathcal{B}(V_{\mathbb{H}}^{R})$. Then, $\sigma_{S}(T+A)\subset\sigma_{S}(A)$ for all $A\in \mathcal{B}(V_{\mathbb{H}}^{R})$ if and only if $T\equiv0$.
\end{theorem}
\proof Assume that $T$ is non-zero operator and $\sigma_{S}(T+A)\subset\sigma_{S}(A)$ for all $A\in \mathcal{B}(V_{\mathbb{H}}^{R})$. Let $0_{V_{\mathbb{H}}^{R}}\neq x\in V_{\mathbb{H}}^{R}$ such that
\begin{align*}Tx=y\neq 0_{V_{\mathbb{H}}^{R}}.\end{align*}
{\bf First step:} There exists $f\in V_{\mathbb{H}}^{'}$ such that
\begin{align*}f(x)=1\mbox{ and }f(y) \neq 0.\end{align*}
Indeed, if $x=yq$ for some $q\in\mathbb{H}$, the result follows from Hahn-Banach theorem. Now assume that $x$ and $y$ are linearly independent. Then, there exists a right basis $Z$ of $V_{\mathbb{H}}^{R}$ such that $\{x,y\}\subset Z$. We consider the following map
\begin{align*}f:u=\displaystyle\sum_{z\in Z}zq_{z}\longmapsto q_{x}+q_{y}.\end{align*}
It is clear that $f$ is right linear and $f(x)=1$ and $f(y)\neq 0$.
\vskip 0.2 cm

\noindent {\bf Second step:} Take
\begin{align*}A:=(x-y)\otimes f\end{align*}
the rank one operator on $V_{\mathbb{H}}^{R}$ given by
\begin{align*}u\longmapsto (x-y)f(u).\end{align*}
We have:
\begin{align*}(T+A)^{2}x-2(T+A)x
&=(T+A)x-2x=-x.\end{align*}
So, $1\in\sigma_{S}(A+T)$. Since $\dim V_{\mathbb{H}}^{R}>1$, then $0\in\sigma_{S}(A)$. In the sequel, we assume that $x\neq y$. Let $q\in\sigma_{S}(A)$. By \cite[Lemma 4.2.3]{DFIS}, we have
\begin{align*}(f(x-y))^{2}-2{\rm Re}(q)f(x-y)+|q|^{2}=0\end{align*}
if and only if
\begin{align*}q\in\Big\{hf(x-y)h^{-1}:\ h\in \mathbb{H}\backslash\{0\}\Big\}.\end{align*}
This implies that $1\not\in\sigma_{S}(A)$. Indeed, if $1$ is an $S-$eigenvalue of $A$, there exists $h\in\mathbb{H}^{*}$ such that
\begin{align*}hf(x-y) h^{-1}=1\end{align*}
\noindent and so $f(y) =0$, contradiction.\qed
\begin{definition}
{\rm A bounded operator $S\in \mathcal{B}(V_{\mathbb{H}}^{R})$ is called a \emph{quasi-inverse} of the operator $T\in \mathcal{B}(V_{\mathbb{H}}^{R})$ if there exists $K_{1},K_{2}\in\mathcal{K}(V_{\mathbb{H}}^{R})$ such that
\begin{align*}ST=\mathbb{I}_{V_{\mathbb{H}}^{R}}-K_{1}\mbox{ and }TS=\mathbb{I}_{V_{\mathbb{H}}^{R}}-K_{2}.\end{align*}}
\end{definition}
\begin{lemma}\label{l:1}
Let $T\in \mathcal{B}(V_{\mathbb{H}}^{R})$, $q\in \mathbb{H}\backslash\sigma_{e}^{S}(T)$. Let $R_{q}(T)$ be a quasi-inverse of $Q_{q}(T)$ and $A$ be an operator that commute with $T$. Then, there exist $K\in \mathcal{K}(V_{\mathbb{H}}^{R})$ such that
\begin{align*}AR_{q}(T)=R_{q}(T)A+K.\end{align*}
\end{lemma}
\proof Since $AT=TA$, then $Q_{q}(T)A=AQ_{q}(T)$. Let $K_{1},\ K_{2}\in \mathcal{K}(V_{\mathbb{H}}^{R})$ such that
\begin{align*}R_{q}(T)Q_{q}(T)=\mathbb{I}_{V_{\mathbb{H}}^{R}}-K_{1}\mbox{  and  }Q_{q}(T)R_{q}(T)=\mathbb{I}_{V_{\mathbb{H}}^{R}}-K_{2}.\end{align*}
Then,
\begin{align*}AR_{q}(T)=R_{q}(T)A+K\end{align*}
where $K=K_{1}AR_{q}(T)-R_{q}(T)AK_{2}\in \mathcal{K}(V_{\mathbb{H}}^{R})$.\qed
\begin{theorem}\label{t:0}
Let $V_{\mathbb{H}}^{R}$ be a  quaternionic Hilbert space and $T\in \mathcal{B}(V_{\mathbb{H}}^{R})$. Let $P_{1}$ be a projector in $\mathcal{B}(V_{\mathbb{H}}^{R})$ commuting with $T$ and let $P_{2}=\mathbb{I}_{V_{\mathbb{H}}^{R}}-P_{1}$. Take $T_{j}:=TP_{j}=P_{j}T,\ j=1,2$. Then,
\begin{align*}\sigma_{e}^{S}(T)=\sigma_{e}^{S}(T_{1}|_{R(P_{1})})\cup\sigma_{e}^{S}(T_{2}|_{R(P_{2})}),\end{align*}
where $T_{i}|_{M}$ denotes the restriction of $T_{i}$ to $M$.
\end{theorem}
\proof Let $q\not\in \sigma_{e}^{S}(T)$. Then, there exist $A_{q}(T)\in \mathcal{B}(V_{\mathbb{H}}^{R})$ and  $K_{1},K_{2}\in \mathcal{K}(V_{\mathbb{H}}^{R})$ such that
\begin{align*}A_{q}(T)Q_{q}(T)=\mathbb{I}_{V_{\mathbb{H}}^{R}}-K_{1}\mbox{ and } Q_{q}(T)A_{q}(T)=\mathbb{I}_{V_{\mathbb{H}}^{R}}-K_{2}.\end{align*}
  Using Lemma \ref{l:1}, we infer that there exists $K_{3}\in \mathcal{K}(V_{\mathbb{H}}^{R})$ such that
 \begin{align*}A_{q}(T)P_{1}=P_{1}A_{q}(T)+K_{3}\mbox{ and }A_{q}(T)P_{2}=P_{2}A_{q}(T)-K_{3}.\end{align*}
  Therefore,
  \begin{equation} \label{e:1}A_{q}(T)=P_{1}A_{q}(T)P_{1}+ P_{2}A_{q}(T)P_{2}+K_{4},\end{equation}
  where $K_{4}=(\mathbb{I}_{V_{\mathbb{H}}^{R}}-2P_{1})K_{3}\in \mathcal{K}(V_{\mathbb{H}}^{R})$. Now, we consider the following identity
\begin{equation}\label{e:2}Q_{q}(T)=(T_{1}^{2}-2{\rm Re}(q)T_{1}+\vert q\vert^{2} P_{1})+(T_{2}^{2}-2{\rm Re}(q)T_{2}+\vert q\vert^{2} P_{2}).\end{equation}
We multiply the identity $(\ref{e:2})$ by $A_{q}(T)$ on the left and on the right, we obtain
\begin{align*}\mathbb{I}_{V_{\mathbb{H}}^{R}}-K_{1}
&=(P_{1}A_{q}(T)P_{1}+ P_{2}A_{q}(T)P_{2}+K_{4})((T_{1}^{2}-2{\rm Re}(q)T_{1}+\vert q\vert^{2} P_{1}))\\
&+(P_{1}A_{q}(T)P_{1}+ P_{2}A_{q}(T)P_{2}+K_{4})(T_{2}^{2}-2{\rm Re}(q)T_{2}+\vert q\vert^{2} P_{2})]\end{align*}
and
\begin{align*}\mathbb{I}_{V_{\mathbb{H}}^{R}}-K_{2}
&=((T_{1}^{2}-2{\rm Re}(q)T_{1}+\vert q\vert^{2} P_{1}))(P_{1}A_{q}(T)P_{1}+P_{2}A_{q}(T)P_{2}+K_{4})\\
&+((T_{2}^{2}-2{\rm Re}(q)T_{2}+\vert q\vert^{2}P_{2})(P_{1}A_{q}(T)P_{1}+P_{2}A_{q}(T)P_{2}+K_{4}) .\end{align*}
 This leads us to conclude that

\begin{align*}\mathbb{I}_{V_{\mathbb{H}}^{R}}-K_{1}
&=P_{1}A_{q}(T)P_{1}(T_{1}^{2}-2{\rm Re}(q)T_{1}+\vert q\vert^{2} P_{1})\\
&+P_{2}A_{q}(T)P_{2}(T_{2}^{2}-2{\rm Re}(q)T_{2}+\vert q\vert^{2} P_{2})+K_{5}\end{align*}
and
\begin{align*}\mathbb{I}_{V_{\mathbb{H}}^{R}}-K_{2}
&=(T_{1}^{2}-2{\rm Re}(q)T_{1}+\vert q\vert^{2} P_{1})P_{1}A_{q}(T)P_{1}\\
&+(T_{2}^{2}-2{\rm Re}(q)T_{2}+\vert q\vert^{2} P_{2})P_{2}A_{q}(T)P_{2}+K_{6},\end{align*}
where
\begin{align*}K_{5}=K_{4}Q_{q}(T)\in \mathcal{K}(V_{\mathbb{H}}^{R})\mbox{ and }K_{6}=Q_{q}(T)K_{4}\in\mathcal{K}(V_{\mathbb{H}}^{R}).\end{align*}

\noindent Take
\begin{align*}A_{q,j}(T)=P_{j}A_{q}(T)P_{j},\mbox{ for }j=1,2.\end{align*}

\noindent Then,
\begin{align*}
P_{1}(\mathbb{I}_{V_{\mathbb{H}}^{R}}-K_{1}-K_{5})P_{1}=A_{q,1}(T)(T_{1}^{2}-2{\rm Re}(q)T_{1}+\vert q\vert^{2} P_{1})
\end{align*}
and
\begin{align*}
P_{1}(\mathbb{I}_{V_{\mathbb{H}}^{R}}-K_{2}-K_{6})P_{1}=(T_{1}^{2}-2{\rm Re}(q)T_{1}+\vert q\vert^{2} P_{1})A_{q,1}(T).
\end{align*}
As a consequence, $A_{q,1}(T)|_{R(P_{1})}$ is a quasi-inverse of $(T_{1}^{2}-2{\rm Re}(q)T_{1}+\vert q\vert^{2} P_{1})|_{R(P_{1})}$. Similarly, we have $A_{q,2}(T)|_{N(P_{1})}$ is a quasi-inverse $(T_{2}^{2}-2{\rm Re}(q)T_{2}+\vert q\vert^{2} P_{2})|_{N(P_{1})}$. Therefore, we conclude that  \begin{align*}q\in\mathbb{H}\backslash(\sigma^{S}_{e}(T_{1}\mid_{R(P_{1})})\cup\sigma^{S}_{e}(T_{2}\mid_{R(P_{2})})).\end{align*}
\vskip 0.1 cm

Conversely, let $q\not\in \sigma_{e}^{S}(T_{1}|_{R(P_{1})})\cup\sigma_{e}^{S}(T_{2}|_{R(P_{2})})$. Then, there exists $A_{q,i}(T)\in\mathcal{B}(R(P_{i}))$ and $K_{i,j}\in\mathcal{K}(R(P_{i}))$, for $i, j=1,2$, such that
\begin{align*}A_{q,i}Q_{q}(T_{i}\mid_{R(P_{i})})=\mathbb{I}_{R(P_{i})}-K_{1,i}\mbox{ and }Q_{q}(T_{i})A_{q,i}(T)=\mathbb{I}_{R(P_{i})}-K_{2,i}.\end{align*}
Let us define the operator
\begin{align*}B_{q}(T):=P_{1}A_{q,1}(T)P_{1}+P_{2}A_{q,2}(T)P_{2}.\end{align*}
\noindent  We get:
\begin{align*}B_{q}(T)Q_{q}(T)
&=P_{1}A_{q,1}(T)Q_{q}(T)P_{1}+P_{2}A_{q,2}(T)Q_{q}(T)P_{2}\\
&=P_{1}(\mathbb{I}_{R(P_{1})}-K_{1,1})P_{1}+P_{2}(\mathbb{I}_{R(P_{2})}-K_{1,2})P_{2}\\
&=P_{1}+P_{2}-P_{1}K_{1,1}P_{1}-P_{2}K_{1,2}P_{2}\\
&=\mathbb{I}_{V_{\mathbb{H}}^{R}}-K_{6}\end{align*}
where $K_{6}=P_{1}K_{1,1}P_{1}+P_{2}K_{1,2}P_{2}\in\mathcal{K}(V_{\mathbb{H}}^{R})$. By a similar argument, we obtain
\begin{align*}Q_{q}(T)B_{q}(T)=\mathbb{I}_{V_{\mathbb{H}}^{R}}-K_{7},\end{align*}
where $K_{7}=P_{1}K_{2,1}P_{1}+P_{2}K_{2,2}P_{2}\in\mathcal{K}(V_{\mathbb{H}}^{R})$. Therefore, $B_{q}(T)$ is a quasi-inverse of $Q_{q}(T)$. We deduce that $q\not\in\sigma_{e}^{S}(T)$.\qed

The result of \cite[Theorem 6]{PSK} can be reformulate as follows:
\begin{theorem}\label{t:1}
Let $T\in\mathcal{B}(V_{\mathbb{H}}^{R})$ and let $\sigma$ be an isolated part of $\sigma_{S}(T)$. Put $V_{1,\mathbb{H}}^{R}=R(P_{\sigma})$ and $V_{2,\mathbb{H}}^{R}=N(P_{\sigma})$. Then,
$V_{\mathbb{H}}^{R}=V_{1,\mathbb{H}}^{R}\oplus V_{2,\mathbb{H}}^{R}$, the spaces $V_{1,\mathbb{H}}^{R}$ and $V_{2,\mathbb{H}}^{R}$ are $T-$invariant subspaces and
\begin{equation}\label{e:3}\sigma_{S}(T|_{V_{1;\mathbb{H}}^{R}})=\sigma\mbox{ and }\sigma_{S}(T|_{V_{2,\mathbb{H}}^{R}})=\sigma_{S}(T)\backslash \sigma.\end{equation}
\end{theorem}
\vskip 0.2 cm

In the next result, we discuss the uniqueness of the decomposition $(\ref{e:3})$.

\begin{theorem}\label{p:1}Let $T\in\mathcal{B}(V_{\mathbb{H}}^{R})$ and $P$ be a projection in $\mathcal{B}(V_{\mathbb{H}}^{R})$ such that $TP=PT$. Set
\begin{align*}T_{1}:=TP\vert_{R(P)}\mbox{ and } T_{2}:=T(\mathbb{I}_{V_{\mathbb{H}}^{R}}-P)\vert_{N(P)}.\end{align*}
\noindent If ${\rm dist}(\sigma_{S}(T_{1}),\sigma_{S}(T_{2}))>0$, then
\begin{align*}R(P)=R(P_{\sigma_{S}(T_{1})})\mbox{ and }N(P)=N(P_{\sigma_{S}(T_{1})})\end{align*}
\end{theorem}

\proof First and foremost, we have by \cite[Theorem 4.4]{KRE}
\begin{align*}\sigma_{S}(T)=\sigma_{S}(T_{1})\cup\sigma_{S}(T_{2})\end{align*}

 \noindent Take $\sigma_{1}=\sigma_{S}(T_{1})$ and $\sigma_{2}=\sigma_{S}(T_{2})$ and assume that ${\rm dist}(\sigma_{1},\sigma_{2})>0$. According to the proof of \cite[Theorem 6]{PSK}, there exist a pair of a non-empty disjoint axially symmetric domains $(U_{\sigma_{1}},U_{\sigma_{2}})$ such that
\begin{align*}\sigma_{j}\subset U_{\sigma_{j}},\ \overline{U}_{\sigma_{1}}\cap\overline{U}_{\sigma_{2}}=\emptyset\end{align*}
and the boundary $\partial(U_{\sigma_{j}}\cap\mathbb{C}_{I})$ is the union of finite number of continuously differentiable Jordan curves for $j=1,2$ and for all $I\in\mathbb{S}$. In this way, we see that $U_{\sigma_{j}}$ is $T_{j}$-admissible open set for $j=1,2$. So by the quaternionic functional calculus, we deduce that
\begin{align*}\displaystyle\frac{1}{2\pi}\int_{\partial(U_{\sigma_{j}}\cap\mathbb{C}_{I})}S_{L}^{-1}(s,T_{j})ds_{I}=\mathbb{I}_{V_{j,\mathbb{H}}^{R}} \mbox{ for }j=1,2\end{align*}
and
\begin{align*}\displaystyle\frac{1}{2\pi}\int_{\partial(U_{\sigma_{i}}\cap\mathbb{C}_{I})}S_{L}^{-1}(s,T_{j})ds_{I}=0 \mbox{ for }i\neq j\end{align*}
where
\begin{align*}V_{1,\mathbb{H}}^{R}:=R(P)\mbox{  and  }V_{2,\mathbb{H}}^{R}:=N(P).\end{align*}
Define the operators:
\begin{align*}R_{1}=\left[
 \begin{array}{ccc}
 \mathbb{I}_{R(P)}& 0 \\
0 & 0 \\
\end{array}
\right] \mbox{ and } R_{2}=\left[
 \begin{array}{ccc}
0& 0 \\
0 &  \mathbb{I}_{N(P)} \\
\end{array}
\right].\end{align*}
We have
\begin{align*}P_{\sigma_{1}}
&=\displaystyle\frac{1}{2\pi}\int_{\partial(U_{\sigma_{1}}\cap\mathbb{C}_{I})}S_{L}^{-1}(s,T)ds_{I}\\
&=\displaystyle\frac{1}{2\pi}\int_{\partial(U_{\sigma_{1}}\cap\mathbb{C}_{I})}ds_{I}S_{R}^{-1}(s,T)\\
&=\displaystyle\frac{1}{2\pi}\int_{\partial(U_{\sigma_{1}}\cap\mathbb{C}_{I})}ds_{I}S_{R}^{-1}(s,T)P+
\displaystyle\frac{1}{2\pi}\int_{\partial(U_{\sigma_{1}}\cap\mathbb{C}_{I})}ds_{I}S_{R}^{-1}(s,T)(\mathbb{I}_{V_{\mathbb{H}}^{R}}-P)\\
&=\displaystyle\frac{1}{2\pi}\int_{\partial(U_{\sigma_{1}}\cap\mathbb{C}_{I})}ds_{I}S_{R}^{-1}(s,T_{1})P+
\displaystyle\frac{1}{2\pi}\int_{\partial(U_{\sigma_{1}}\cap\mathbb{C}_{I})}ds_{I}S_{R}^{-1}(s,T_{2})(\mathbb{I}_{V_{\mathbb{H}}^{R}}-P)\\
&=\displaystyle\frac{1}{2\pi}\int_{\partial(U_{\sigma_{1}}\cap\mathbb{C}_{I})}ds_{I}S_{R}^{-1}(s,T_{1})P\\
&=P_{\sigma_{1}}P\\
&=(I-P_{\sigma_{2}})P\\
&=R_{1}.\end{align*}
We conclude that
\begin{align*}R(P)=R(P_{\sigma_{1}}).\end{align*}

\noindent By the similar arguments, we achieve that
\begin{align*}P_{\sigma_{2}}=R_{2}\end{align*}
 and so, $N(P)=R(P_{\sigma_{2}})=N(P_{\sigma_{1}})$.\qed
\vskip 0.3 cm

We now analyze the Riesz projection associated with the isolated spheres. To start off, we give the following definition:
\begin{definition}
{\rm $T\in\mathcal{B}(V_{\mathbb{H}}^{R})$.  A point $q\in\sigma_{S}(T)$ is called an eigenvalue of finite type if $V_{\mathbb{H}}^{R}$ is a direct sum of $T-$invariant subspaces $V_{1,\mathbb{H}}^{R}$
and $V_{2,\mathbb{H}}^{R}$ such that
\vskip 0.1 cm

$(H1)$ $\dim(V_{1,\mathbb{H}}^{R})<\infty$,
\vskip 0.1 cm

$(H2)$ $\sigma_{S}(T\vert_{V_{1,\mathbb{H}}^{R}})\cap \sigma_{S}(T\vert_{V_{2,\mathbb{H}}^{R}})=\emptyset$,
\vskip 0.1 cm

$(H3)$ $\sigma_{S}(T\vert_{V_{1,\mathbb{H}}^{R}})=[q]$.}

\end{definition}
\begin{remark}~~\\
\begin{enumerate}
{\rm \item In complex spectral theory, in a complex Hilbert space $V_{\mathbb{C}}$, for a continuous linear operator, $T$, the condition $(H3)$ is replaced by $\sigma(T\vert_{V_{1,\mathbb{C}}})=\{q\}$ $($where $V_{\mathbb{C}}$ is a direct sum of $T-$invariant subspaces $V_{1,\mathbb{C}}$ and $V_{2,\mathbb{C}})$, see \cite{Gohberg}. In the quaternionic setting, we must take the whole $2-$sphere $[q]$ because if $q\in \sigma_{S}(T\vert_{V_{1,\mathbb{H}}^{R}})$, then $[q]\subset \sigma_{S}(T\vert_{V_{1,\mathbb{H}}^{R}})$.
\item Let $T\in\mathcal{B}(V_{\mathbb{H}}^{R})$. If $q\in\sigma_{S}(T)\backslash\mathbb{R}$, then $q$ is not an isolated point of $\sigma_{S}(T)$. Take
\begin{align*}\Omega:=\sigma_{S}(T)/\sim\end{align*}
where $p\sim q$ if and only if $p\in [q]$. Let $E_{T}$ denote the set of representatives of $\Omega$.
 \item By \cite[Proposition 4.44 and Theorem 4.47]{KM}, If $\dim(V_{\mathbb{H}}^{R})<\infty$ , then the $S-$spectrum of $T$ consists of right eigenvalues only and $\# E_{T}<\infty$. In particular, if $T$ satisfied the assumptions $(H1)$ and $(H3)$, then  $q\in\sigma_{pS}(T)$.}
\end{enumerate}
\end{remark}
\begin{lemma}
Let $T\in\mathcal{B}(V_{\mathbb{H}}^{R})$. Let $E_{T}$ denote the set of representatives as above. Then, $q$ is an isolated point of $E_{T}$ if and only if $[q]$ is an isolated part of $\sigma_{S}(T)$.
\end{lemma}
\proof Assume that $[q]$ is isolated part of $\sigma_{S}(T)$ and let $U_{q}\subset\mathbb{H}$ be an open set such that
\begin{align*}[q]=\sigma_{S}(T)\cap U_{q}.\end{align*}
If there exists $q\neq p\in E_{T}\cap U_{q}$, then $p\not\in [q]$. This implies that
\begin{align*}E_{T}\cap U_{q}=\{q\}.\end{align*}
Conversely, if $q$ is an isolated point of $E_{T}$, then $[q]$ is an open subset of $\sigma_{S}(T)$. Since $[q]$ is compact, then $[q]$ is isolated part of $\sigma_{S}(T)$.\qed
\vskip 0.3 cm

Let $V_{\mathbb{C}}$ be a complex Hilbert space and $T$ be a continuous linear operator acting on $V_{\mathbb{C}}$. It follows from \cite[Theorem 1.1]{Gohberg} that if $q\in\sigma(T)$ $($where $\sigma(T)$ denotes the complex spectrum of $T)$, then $q$ is an eigenvalue of finite type if and only if $\dim R(P_{\{q\}})<\infty$, where $P_{\{q\}}$ is the Riesz projection corresponding to the isolated point $q$. In the next theorem, we show that the same is true for the right eigenvalue of finite type in the quaternionic setting.
\begin{theorem}\label{t:2}
Let $T\in\mathcal{B}(V_{\mathbb{H}}^{R})$ and $q\in\sigma_{pS}(T)$. Then, $q$ is a right eigenvalue of finite type if and only if $\{q\}$ is an isolated part in $E_{T}$ and $\dim R(P_{[q]})<\infty$.
\end{theorem}
\proof
Assume that $q$ is a right eigenvalue of finite type and consider the direct sum
\begin{align*}V_{\mathbb{H}}^{R}=V_{1,\mathbb{H}}^{R}\oplus V_{2,\mathbb{H}}^{R}\end{align*}
 with the properties $(H1)-(H3)$. By Theorem \ref{p:1}, this decomposition is unique and so
\begin{align*}V_{1,\mathbb{H}}^{R}= R(P_{[q]}).\end{align*}
The converse comes from Theorem \ref{t:1}.
\begin{corollary}
Let $T\in\mathcal{B}(V_{\mathbb{H}}^{R})$. Then,
\begin{align*}\sigma_{e}^{S}(T)\subset \sigma_{S}(T)\backslash\sigma_{d}^{S}(T)\end{align*}
where $\sigma_{d}^{S}(T)$ denotes the set of all right eigenvalue of finite type.
\end{corollary}

\proof Assume that $q$ is a right eigenvalue of finite type. Then, $[q]$ is an isolated part of $\sigma_{S}(T)$ and $\dim R(P_{[q]})<\infty$. Therefore,
\begin{align*}\dim R(TP_{[q]})<\infty\mbox{ and so }\sigma_{e}^{S}(T\vert_{R(P_{[q]})})=\emptyset.\end{align*}
Using Theorem \ref{t:0}, we infer that
\begin{align*}\sigma_{e}^{S}(T)=\sigma_{e}^{S}(T\vert_{R(\mathbb{I}_{V_{\mathbb{H}}^{R}}-P_{[q]})}).\end{align*}
On the other hand,
\begin{align*}\sigma_{e}^{S}(T\vert_{R(\mathbb{I}_{V_{\mathbb{H}}^{R}}-P_{[q]})})\subset \sigma_{S}(T\vert_{R(\mathbb{I}_{V_{\mathbb{H}}^{R}}-P_{[q]})})=\sigma_{S}(T)\backslash [q].\end{align*}
This leads to conclude that $q\not\in\sigma_{e}^{S}(T)$.\qed
\vskip 0.3 cm

\begin{example}
{\rm We consider the right quaternionic Hilbert space:
\begin{align*}\ell^{2}_{\mathbb{H}}(\mathbb{Z}):=\Big\{x:\mathbb{Z}\longrightarrow\mathbb{H}\mbox{  such that  } \|x\|^{2}:=\sum_{i\in\mathbb{Z}}|x_{i}|^{2}<\infty\Big\}.\end{align*}
with the right  scalar multiplication
\begin{align*}xa=(x_{i}a)_{i\in\mathbb{Z}}\end{align*}
for $x=(x_{i})_{i\in\mathbb{Z}}$ and $a\in\mathbb{H}$. The associated scalar product is given by
\begin{align*}\langle x,y\rangle:=\langle x,y\rangle_{\ell^{2}_{\mathbb{H}}(\mathbb{Z})}:=\sum_{i\in\mathbb{Z}}\overline{x_{i}}y_{i}.\end{align*}
The right shift is the map:
\begin{align*}
T:\ &\ell_{\mathbb{H}}^{2}(\mathbb{Z})\longrightarrow \ell_{\mathbb{H}}^{2}(\mathbb{Z})
\\
& x\longmapsto\ y=(y_{i})_{i\in\mathbb{Z}}
\end{align*}
where $y_{i}=x_{i+1}$ if $i\neq -1$ and $0$ if $i=-1$. We have
\begin{align*}\|T(x)\|^{2}=\sum_{i\neq -1}|x_{i}|^{2}\leq \|x\|^{2}.\end{align*}
 The $S-$spectrum of $T$ was studied in \cite{B1,BK4}. In particular, we have
\begin{align*}\sigma_{S}(T)=\sigma_{pS}(T)=\nabla_{\mathbb{H}}(0,1),\end{align*}
where $\nabla_{\mathbb{H}}(0,1)$ is the closed quaternionic unit ball. By Theorem \ref{t:2}, none of these $S-$eigenvalues are of finite type.}
\end{example}

We recall:
\begin{lemma}\label{l:2}\cite[Corollary 2.22]{GN0}
Let $V_{\mathbb{H}}^{R}$ be a  quaternionic right vector space and $E$ be a right linear independent subspace of $V_{\mathbb{H}}^{R}$. Then, there exists a right basis $B$ of $V_{\mathbb{H}}^{R}$ such that $E\subset B$. In particular, every quaternionic right vector space has a right basis.
\end{lemma}
\begin{proposition}\label{p:2}\cite[Proposition 4.44]{KM}
Let $V_{\mathbb{H}}^{R}$ be a quaternionic Hilbert space and $T\in\mathcal{B}(V_{\mathbb{H}}^{R})$. If $q_{1},....,q_{n}\in\mathbb{H}$ are right eigenvalues of $T$ such that $[q_{i}]\neq[q_{j}]$, $\forall 1\leq i<j\leq n\ (n\geq 2)$, and $Q_{q_{j}}(T)x_{j}=0,\ 0\neq x_{j}\in V_{\mathbb{H}}^{R},\forall 1\leq j \leq n$, then $x_{1},x_{2},...,x_{n}$ are right-linearly independent in $V_{\mathbb{H}}^{R}$.
\end{proposition}
\vskip 0.1 cm

We have the following lemma:
\begin{lemma}\label{l:3}
Assume that $\dim(V_{\mathbb{H}}^{R})<\infty$ and let  $T\in\mathcal{B}(V_{\mathbb{H}}^{R})$. Then, $\# E_{T}<\infty$. In this case, set \begin{align*}E_{T}=\Big\{q_{1},q_{2},...,q_{n}\Big\},\end{align*}
then $V_{\mathbb{H}}^{R}$ is a direct sum of $T-$invariant right subspaces $V_{1,\mathbb{H}}^{R},\ V_{2,\mathbb{H}}^{R},\  ....,\ V_{n,\mathbb{H}}^{R}$. Moreover, if $T_{i}:=T\vert_{V_{i,\mathbb{H}}^{R}}:\ V_{i,\mathbb{H}}^{R}\longrightarrow V_{i,\mathbb{H}}^{R}$, for $i=1,...,n$, then

\begin{align*}\sigma_{S}(T_{i})=\Big\{hq_{i}h^{-1}:\ h\in\mathbb{H}^{*}\Big\}.\end{align*}
\end{lemma}
\proof Combine Lemma \ref{l:2}, Proposition \ref{p:2} and Theorem \ref{t:2}\qed
\vskip 0.1 cm

Let $V_{\mathbb{H}}^{R}$ be a quaternionic right Hilbert space, $T\in\mathcal{B}(V_{\mathbb{H}}^{R})$ and $q$ be a right eigenvalue of $T$ of finite type. In this way, we see that
\begin{align*}V_{\mathbb{H}}^{R}=R(P_{[q]})\oplus R(P_{\sigma_{S}(T)\backslash [q]}).\end{align*}

\begin{definition}{\rm The algebraic multiplicity of the right eigenvalue $q$ is, by definition, the dimension of the space $R(P_{[q]})$.}\end{definition}
In the next, we write:
\begin{align*}m_{T}(q):=\dim R(P_{[q]}).\end{align*}
We refer to \cite{Gohberg,J1} for the definition in complex setting. Now, inspired by the description of the Riesz projection in the complex setting, see \cite{Gohberg}, we give a quaternionic version of a result describing the finite part of right eigenvalues of finite type.
\begin{theorem}Let $V_{\mathbb{H}}^{R}$ be a quaternionic Hilbert space, $T\in\mathcal{B}(V_{\mathbb{H}}^{R})$ and $\sigma$
be an axially symmetric isolated part of $\sigma_{S}(T)$. Set $P_{\sigma}$ the Riesz projection correspond to $\sigma$ and take
\begin{align*}E_{T}^{\sigma}=\sigma/\cong\end{align*}
where $q\cong p$ if and only if $p\in [q]$. Then, $\dim R(P_{\sigma})<\infty$ if and only if $\# E_{T}^{\sigma}<\infty$ and $q$ is a right eigenvalue of finite type for all $q\in E_{T}^{\sigma}$. Besides, if so, then
\begin{align*}\dim R(P_{\sigma})=\displaystyle\sum_{q\in E_{T}^{\sigma}}\dim R(P_{[q]}).\end{align*}
\end{theorem}
\proof
If $\dim R(P_{\sigma})<\infty$, then $P_{\sigma}T$ is a finite rank operator. In this way, we see that
\begin{align*}\# E_{T|_{R(P_{\sigma})}}^{\sigma}=n<\infty.\end{align*}
By Lemma \ref{l:3},
\begin{align*}R(P_{\sigma})=V_{1,\mathbb{H}}^{R}\oplus V_{2,\mathbb{H}}^{R}\oplus ...\oplus V_{n,\mathbb{H}}^{R}\end{align*}
where $V_{j,\mathbb{H}}^{R}$ is $T-$invariant with the properties
\begin{align*}\sigma_{S}(T|_{V_{j,\mathbb{H}}^{R}})=[q_{j}],\ j=1,...,n.\end{align*}
Let $i\in\{1,...,n\}$. Since $P_{\sigma}$ is a projection, then
\begin{align*}V_{\mathbb{H}}^{R}=N( P_{\sigma})\oplus R(P_{\sigma}).\end{align*}
In this fashion, we have
\begin{align*}V_{\mathbb{H}}^{R}=V_{i,\mathbb{H}}^{R}\oplus W_{i,\mathbb{H}}^{R}\end{align*}
where
\begin{align*}W_{i,\mathbb{H}}^{R}=V_{1,\mathbb{H}}^{R}\oplus ...\oplus V_{i-1,\mathbb{H}}^{R}\oplus V_{i+1,\mathbb{H}}^{R}\oplus...\oplus V_{n, \mathbb{H}}^{R}\oplus  N(P_{\sigma}).\end{align*}

As consequence, $V_{i,\mathbb{H}}^{R}$ and $W_{i,\mathbb{H}}^{R}$ are $T-$invariant and $\sigma_{S}(T|_{V_{i,\mathbb{H}}^{R}})=[q_{i}]$. So, $q_{i}$ is a right eigenvalue of finite type.
\skip 0.1 cm

Conversely, set $E_{T}^{\sigma}=\{q_{1},\ q_{2}...,q_{n}\}$, where $q_{i}$ is a right eigenvalue of $T$ of finite type for all $i\in\{1,2,...,n\}$.
 Applying \cite[Theorem 5.6]{NCFCBO}, we have
 \begin{align*}\mathbb{I}_{R(P_{\sigma})}=\displaystyle \sum_{i=1}^{n}P_{[q_{i}]}|_{R(P_{\sigma})}.\end{align*}
Since $P_{[q_{i}]}P_{[q_{j}]}=0$ for all $i\neq j$, then
\begin{align*}R(P_{\sigma})=R(P_{[q_{1}]})\oplus R(P_{[q_{2}]})\oplus ...\oplus R(P_{[q_{n}]}).\end{align*}
This implies that
\begin{align*}\dim (R(P_{\sigma}))=\displaystyle\sum_{i=1}^{n}\dim (R(P_{[q_{i}]})).\end{align*}
In particular, we have $P_{\sigma}$ which is a finite rank operator.\qed
\begin{remark}
{\rm In the complex spectral theory, much attention has ben paid to eigenvalue of finite type, see \cite{charfi,Gohberg,J1,J2,Lutgen}. It is useful for the study of the essential spectrum of certain operators-matrices. We refer to \cite{charfi} for this point on the two-groupe transport operators. More precisely, let $V_{\mathbb{C}}$ be a complex Banach space and let $T$ be a closed operator in $V_{\mathbb{C}}$. The Browder resolvent set of $T$ is given by
\begin{align*}\rho_{B}(T):=\rho(T)\cup\sigma_{d}(T),\end{align*}
where we use the notation $\rho(.)$ for the resolvent set of $T$ and $\sigma_{d}(.)$ the set of eigenvalues of finite type of $T$. In fact, the usual resolvent
\begin{align*}R_{\lambda}(A):=(A-\lambda)^{-1}\end{align*}
can be extended to $\rho_{B}(T)$, e.g. \cite{Lutgen}. Motivated by this,  \cite{charfi} gives a version of the Frobenius-Schur factorization  using the Browder resolvent. This makes it possible to study the essential spectrum of serval types operators-matrices. In this paper, we have described the discrete $S-$spectrum in quaternionic setting. In this regard, as in complex case, we can define the spherical Browder resolvent. Although, we avoided studying it in this paper, we will cover that in a future article.}
\end{remark}
\section{Some results on the Weyl $S$-spectrum}\label{sec:3}
In this section, we develop a deeper understanding of the concept of the Weyl $S$-spectrum of the bounded right linear operator. More precisely, we describe the boundary of the $S$-spectrum. Likewise, we deal with the particular case of the spectral theorem. To begin with, we recall:
\begin{definition}\cite{BK2}{\rm Let $T\in\Bc(V_{\mathbb{H}}^{R})$. The Weyl $S-$spectrum is the set
\begin{align*}\sigma_{W}^{S}(T)=\displaystyle \bigcap_{K\in\mathcal{K}(V_{\mathbb{H}}^{R})}\sigma_{S}(T+K).\end{align*}}
\end{definition}
\vskip 0.1 cm

The study of the essential and the Weyl $S-$spectra are established using the Fredholm theory, see \cite{BK,BK2}. We refer to \cite{B1} for the investigation of the Fredholm and Weyl elements with respect to a quaternionic Banach algebra homomorphism.
\begin{definition}
{\rm A Fredholm operator is an operator $T\in\Bc(V_{\mathbb{H}}^{R})$ such that $N(T)$ and $V_{\mathbb{H}}^{R}/R(T)$ are finite dimensional. We will denote by $\Phi(V_{\mathbb{H}}^{R})$ the set of all Fredholm operators.}
\end{definition}
\vskip 0.1 cm

From \cite{BK,BK2}, we have
\begin{align*}\Phi(V_{\mathbb{H}}^{R})=\Phi_{l}(V_{\mathbb{H}}^{R})\cap\Phi_{r}(V_{\mathbb{H}}^{R})\end{align*}
where
\begin{align*}\Phi_{l}(V_{\mathbb{H}}^{R})=\Big\{T\in \Bc(V_{\mathbb{H}}^{R}):\mbox{ R(T) is closed and }\dim (N(T))<\infty\Big\}\end{align*}
and
\begin{align*}\Phi_{r}(V_{\mathbb{H}}^{R})=\Big\{T\in \Bc(V_{\mathbb{H}}^{R}):\mbox{ R(T) is closed and }\dim (N(T^{\dag}))<\infty\Big\}.\end{align*}
Let $T\in \Phi_{l}(V_{\mathbb{H}}^{R})\cup \Phi_{r}(V_{\mathbb{H}}^{R})$. Then, the index of $T$ is given by
\begin{align*}i(T):=\dim N(T)-\dim(V_{\mathbb{H}}^{R}/R(T)).\end{align*}
\begin{theorem}\cite{BK,BK2} \label{t:5}
Let $T\in\Bc(V_{\mathbb{H}}^{R})$. Then,
\begin{align*}\sigma_{e}^{S}(T)=\mathbb{H}\backslash \Phi_{T}\mbox{  and  }\sigma_{W}^{S}(T)=\mathbb{H}\backslash W_{T}\end{align*}
where
\begin{align*}\Phi_{T}:=\Big\{q\in\mathbb{H}:\ Q_{q}(T)\in\Phi(V_{\mathbb{H}}^{R})\Big\}\end{align*}
\mbox{ and }
\begin{align*}W_{T}:=\Big\{q\in\mathbb{H}:\ Q_{q}(T)\in\Phi(V_{\mathbb{H}}^{R})\mbox{ and }i(Q_{q}(T))=0\Big\}.\end{align*}
\end{theorem}
\begin{remark} {\rm Let $V_{\mathbb{H}}^{R}$ be a quaternionic space and $T\in\mathcal{B}(V_{\mathbb{H}}^{R})$.
\begin{enumerate}
\item Note that, in general, we have
\begin{align*}\sigma_{e}^{S}(T)\subset\sigma_{W}^{S}(T)=\sigma_{1,W}^{S}(T)\cup\sigma_{2,W}^{S}(T)\subset \sigma_{S}(T)\backslash\sigma_{d}(T).\end{align*}
where
\begin{align*}\sigma_{1,W}^{S}(T):=\mathbb{H}\backslash\Big\{q\in\mathbb{H}:\ Q_{q}(T)\in \Phi_{l}(V_{\mathbb{H}}^{R})\mbox{ and }i(Q_{q}(T))\leq 0\Big\}\end{align*}
and
\begin{align*}\sigma_{2,W}^{S}(T):=\mathbb{H}\backslash\Big\{q\in\mathbb{H}:\ Q_{q}(T)\in \Phi_{r}(V_{\mathbb{H}}^{R})\mbox{ and }i(Q_{q}(T))\geq 0\Big\}.\end{align*}
\noindent In particular, $\sigma_{W}^{S}(T)$ does not contain eigenvalues of finite type.
\item In \cite{B1}, one proves that $q\longmapsto i(T)$ is constant on any component of $\Phi_{T}$. In this way, we see that if $\Phi_{T}$ is connected, then
\begin{align*}\sigma_{e}^{S}(T)=\sigma_{W}^{S}(T).\end{align*}
\end{enumerate}}
\end{remark}
\vskip 0.3 cm

\noindent The first result in this section is the next theorem.
\begin{theorem}\label{t:4}
Let $T\in \mathcal{B}(V_{\mathbb{H}}^{R})$. Then,
\begin{align*}\partial\sigma_{W}^{S}(T)\subset\sigma_{1,W}^{S}(T).\end{align*}
In particular, if $\Phi_{T}$ is connected, then
\begin{align*}\partial\sigma_{e}^{S}(T)=\partial\sigma_{W}^{S}(T)\subset\Big\{q\in\mathbb{H}:\ Q_{q}(T)\not\in \Phi_{l}(V_{\mathbb{H}}^{R})\Big\}.\end{align*}
\end{theorem}
To prove Theorem \ref{t:4}, we first study the concept of the minimum modulus. Let $V_{\mathbb{H}}^{R}$ be a separable right Hilbert space and $T\in \Bc(V_{\mathbb{H}}^{R})$. The minimum modulus of $T$ is given by
 \begin{align*}\mu(T):=\displaystyle\inf_{\|x\|=1}\|Tx\|.\end{align*}
 \vskip 0.3 cm

 \noindent To begin with, we give the following lemma.
 \begin{lemma}\label{l1}
 Let $T$ and $S$ be two bounded right linear operators on a right quaternionic Hilbert space. Then,

 \begin{enumerate}
 \item  If $\|T-S\|<\mu(T)$, then $\mu(S)>0$ and $\overline{R(S)}$ is not a proper subset of $\overline{R(T)}$.\\
 \item  If $\|T-S\|<\frac{\mu(T)}{2}$, then $\overline{R(S)}$ is not a proper subset of $\overline{R(T)}$ and $\overline{R(T)}$
  is not a proper subset of $\overline{R(S)}$.

 \end{enumerate}
 \end{lemma}
 \proof The proof is the same as for the complex Banach space, see Lemma 2.3 and lemma 2.4 in \cite{HAE} for a complex proof.\qed
 \vskip 0.3 cm

For $T\in \mathcal{B}(V_{\mathbb{H}}^{R})$, $q\in\mathbb{H}$ and $\varepsilon>0$ we set:

 \begin{align*}\Oc(T,q,\varepsilon):=\Big\{q'\in\mathbb{H}:\ 2\vert {\rm Re}(q)-{\rm Re}(q')\vert\|T\|+\vert |q'|^{2}-|q|^{2}\vert<\varepsilon.\Big\}\end{align*}
 It is clear  that $\Oc(T,q,\varepsilon)$ is an open set in $\mathbb{H}$.
 \begin{corollary}Let  $T\in \Bc(V_{\mathbb{H}}^{R})$ and $q_{0}\in \rho_{S}(T)$. Then, $q\in \rho_{S}(T)$ for each $q\in \Oc(T,q_{0},\mu(Q_{q_{0}}(T)))$.
\end{corollary}
\proof Let $q\in \Oc(T,q_{0},\mu(Q_{q_{0}}(T)))$. Then,

\begin{align*}\|Q_{q}(T)-Q_{q_{0}}(T)\|
&=\|2({\rm Re}(q_{0})-{\rm Re}(q))T+|q|^{2}-|q_{0}|^{2}\|\\
&\leq 2\vert ({\rm Re}(q_{0})-{\rm Re}(q)\vert\|T\|+\vert|q|^{2}-|q_{0}|^{2} \vert\\
&<\mu(Q_{q_{0}}(T)).\end{align*}

We can apply Lemma \ref{l1} to conclude that

\begin{align*}\mu(Q_{q}(T))>0\mbox{ and }\overline{R(Q_{q}(T))}=\overline{R(Q_{q_{0}}(T))}=V_{\mathbb{H}}^{R}.\end{align*}

By \cite[Proposition 3.5]{BK}, $R(Q_{q}(T))$ is closed. Hence, $q\in \rho_{S}(T)$.\qed

\noindent We recall:
\begin{lemma}\cite[Lemma 7.3.9]{DFIS} \label{l2}Let $n\in\mathbb{N}$ and $q,s\in\mathbb{H}$. Set
\begin{align*}P_{2n}(q)=q^{2n}-2{\rm Re}(s^{n})q^{n}+|s^{n}|^{2}.\end{align*}
Then,
\begin{align*}P_{2n}(q)
&=\Qc_{2n-2}(q)(q^{2}-2{\rm Re}(s)q+|s|^{2})\\
&=(q^{2}-2{\rm Re}(s)q+|s|^{2})\Qc_{2n-2}(q),\end{align*}
where $\Qc_{2n-2}(q)$ is a polynomial of degree $2n-2$ in $q$.
\end{lemma}
\vskip 0.3 cm

\noindent \emph{\bf Proof of Theorem \ref{t:4}}
\noindent Set:
\begin{align*}f(T):=\displaystyle\sup_{K\in \mathcal{K}(V_{\mathbb{H}}^{R})}\mu(T+K).\end{align*}
 Similar proof in the complex case, we have $f(T)>0$ if and only if
 \begin{align*}T\in\Phi_{l}(V_{\mathbb{H}}^{R})\mbox{ and }\dim N(T)\leq\dim(V_{\mathbb{H}}^{R}/R(T)).\end{align*}

  Now,
   Since $\sigma_{e}^{S}(T)$ is not empty (e.g., \cite[Proposition 7.14]{BK}) and $\sigma_{e}^{S}(T)\subset\sigma_{W}^{S}(T)$, then $\partial\sigma_{W}^{S}(T)$ is not empty. Let us then take an element $p$ in $\partial\sigma_{W}^{S}(T)$. Assume that $p\not \in\sigma_{1,w}^{S}(T)$. Then,
   \begin{align*} f(Q_{p}(T))>0.\end{align*}
    So, there exists $K_{0}\in\mathcal{K}(V_{\mathbb{H}}^{R})$ such that
   \begin{align*}\mu(Q_{p}(T)+K_{0})>0.\end{align*}

 \noindent Since $\Oc(T,p,\frac{\mu(Q_{p}(T)+K_{0})}{2})$ is an open neighborhood of $p$ and
\begin{align*}p\in\overline{\Big\{q\in\mathbb{H}:\ Q_{q}(T)\in\Phi(V_{\mathbb{H}}^{R})\mbox{ and }i(Q_{q}(T))=0\Big\}},\end{align*}
then there exists $p_{0}\in \Oc(T,p,\frac{\mu(Q_{p}(T)+K_{0})}{2})$ such that
\begin{align*}p_{0}\in W_{T}.\end{align*}
\noindent On the other hand,
\begin{align*}\|Q_{p_{0}}(T)+K_{0}-Q_{p}(T)-K_{0}\|
&\leq \vert |p_{0}|^{2}-|p|^{2}\vert+2\vert Re(p)-Re(p_{0})\vert\|T\|\\
&<\displaystyle \frac{\mu(Q_{p}(T)+K_{0})}{2}.\end{align*}
Applying Lemma \ref{l1}, we obtain
\begin{align*}R(Q_{p}(T)+K_{0})=V_{\mathbb{H}}^{R} .\end{align*}
Indeed, since $\mu(Q_{p_{0}}+K_{0})>0$ and $p_{0}\in W_{T}$, then
\begin{align*}\dim(V_{\mathbb{H}}^{R}/R(Q_{p_{0}}(T)+K_{0}))=0.\end{align*}
In this way, we see that
\begin{align*}Q_{p}(T)+K_{0}\in\Phi(V_{\mathbb{H}}^{R})\mbox{ and }i(Q_{p}(T)+K_{0})=0.\end{align*}
This implies that, $p\not\in \sigma_{W}^{S}(T)$.
\vskip 0.1 cm
\noindent The rest of the proof follows immediately from \cite[Theorem 5.13]{B1}.\qed
\vskip 0.3 cm

We will now deal with the particular spectral theorem for the essential  S-spectra.

\begin{theorem}
Let $T\in\Bc(V_{\mathbb{H}}^{R})$. Then,
\begin{align*}\sigma_{e}^{S}(T^{n})=\Big\{q^{n}\in\mathbb{H}:\ q\in\sigma_{e}^{S}(T)\Big\}=(\sigma_{e}^{S}(T))^{n}.\end{align*}
\end{theorem}
\proof
According to \cite[Lemma 3.10]{DFIS} and the proof of \cite[Theorem 7.3.11]{DFIS} we have
\begin{align*}T^{2n}-2{\rm Re}(q)T^{n}+|q|^{2}\mathbb{I}_{V_{\mathbb{H}}^{R}}=\displaystyle\prod_{j=0}^{n-1}
(T^{2}-2Re(q_{j})T+|q_{j}|^{2}\mathbb{I}_{V_{\mathbb{H}}^{R}}),\end{align*}
where $q_{j}, j=0,...,n-1$ are the solutions of $p^{n}=q$ in the complex plane $\mathbb{C}_{I_{q}}$. Let $q\in\sigma_{e}^{S}(T^{n})$. Then, $Q_{q}(T^{n})\not\in\Phi(V_{\mathbb{H}}^{R})$.
 We can apply \cite[Theorem 6.13]{BK}, we infer that there exists $i\in\{0,1,...,n-1\}$ such that
\begin{align*}Q_{q_{i}}(T)\not\in \Phi(V_{\mathbb{H}}^{R}).\end{align*}
Therefore, $q_{i}\in\sigma_{e}^{S}(T)$. In this way, we see that $q=q_{i}^{n}\in(\sigma_{e}^{S}(T))^{n}$. To prove the inverse inclusion, we consider $p=q^{n}$, where $q\in\sigma_{e}^{S}(T)$.
By Lemma \ref{l2} and \cite[Theorem 7.3.7]{DFIS}, we get
\begin{align*}T^{2n}-2{\rm Re}(q^{n})T^{n}+\vert q^{n}\vert^{2}\mathbb{I}_{V_{\mathbb{H}}^{R}}
&=\Qc_{2n-2}(T)(T^{2}-2{\rm Re}(q)T+|q|^{2}\mathbb{I}_{V_{\mathbb{H}}^{R}})\\
&=(T^{2}-2{\rm Re}(q)T+|q|^{2}\mathbb{I}_{V_{\mathbb{H}}^{R}})\Qc_{2n-2}(T).\end{align*}
\noindent Since $Q_{q}(T)\not\in \Phi(V_{\mathbb{H}}^{R})$, we can apply \cite[Corollary 6.14]{BK}, we deduce that
\begin{align*}T^{2n}-2{\rm Re}(q^{n})T^{n}+\vert q^{n}\vert^{2}\mathbb{I}_{V_{\mathbb{H}}^{R}}\not\in \Phi(V_{\mathbb{H}}^{R}).\end{align*}
So, $p\in\sigma_{e}^{S}(T^{n})$. \qed

\end{document}